\documentclass[12pt]{article}
\usepackage[a4paper, total={7in, 10in}]{geometry}
\usepackage[parfill]{parskip}
\usepackage{physics, tensor, float, subcaption}
\usepackage{graphicx}
\graphicspath{ {Plots/} }
\usepackage{jhep-mod}
\usepackage{bm}
\usepackage{soul}
\usepackage{amssymb,amsmath,amsthm}
\usepackage{mathrsfs}
\usepackage[utf8]{inputenc}
\usepackage{enumerate}
\usepackage{bigints}
\usepackage{xcolor}
\usepackage{appendix}
\usepackage{graphicx}
\usepackage{float}
\usepackage{tikz}
\usepackage{setspace}
\usepackage{cancel}
\usepackage{doi}
\definecolor{purple}{rgb}{1,0,1}
\definecolor{lime}{HTML}{A6CE39} 

\newcommand{\blue}[1]{{\color{blue} #1}}

\definecolor{lime}{HTML}{A6CE39}
\newcommand{\orcidicon}{%
	\begin{tikzpicture}
	\draw[lime, fill=lime] (0,0) 
		circle [radius=0.16] 
		node[white] {{\fontfamily{qag}\selectfont \tiny ID}};
	\draw[white, fill=white] (-0.0625,0.095) 
		circle [radius=0.007];
	\end{tikzpicture}
	\hspace{-5mm}
}
\newcommand\orcidMatt{{\href{https://orcid.org/0000-0003-1088-6485}{\orcidicon}}}

\renewcommand{\O}{\mathcal{O}}
\def\theta{{\vartheta}}
\def\half{{\textstyle{1\over2}}}

\def\quarter{{\textstyle{1\over4}}}

\def\Heaviside{{\mathrm{Heaviside}}}
\begin{document}

\title{
{\blue{Behaviour of the sequence 
$\theta_n = \theta(p_n)$}}
}

\author{
\Large
Matt Visser\!\orcidMatt\!
}
\affiliation{School of Mathematics and Statistics, Victoria University of Wellington, \\
\null\qquad PO Box 600, Wellington 6140, New Zealand.}
\emailAdd{matt.visser@sms.vuw.ac.nz}
\def\theta{\vartheta}

\abstract{
\vspace{1em}

The well-known sequence $\theta_n = \theta(p_n) = \sum_{i=1}^n \ln p_i= \ln\left([p_n]\#\right)$ exhibits numerous extremely interesting properties. 
Since $p_n = \exp(\theta_n - \theta_{n-1})$, it is immediately clear that the two sequences $p_n \longleftrightarrow \theta_n$ must ultimately encode exactly the same information. But the sequence $\theta_n$,  while being extremely closely correlated with the primes, (in fact, $\theta_n \sim p_n$),  is very much better behaved than the primes themselves. \break
Using numerous suitable extensions of various reasonably standard results,  I shall demonstrate that the sequence $\theta_n$ satisfies suitably defined $\theta$-analogues of the usual Cramer, Andrica, Legendre, Oppermann, Brocard, Firoozbakht, Fourges, Nicholson, and Farhadian  conjectures. (So these $\theta$-analogues are not conjectures, they are instead theorems.) The crucial key to enabling this pleasant behaviour is the regularity (and relative smallness) of  the $\theta$-gaps $\mathfrak{g}_n = \theta_{n+1}-\theta_n= \ln p_{n+1}$. 
While superficially these results bear close resemblance to some recently derived results for the averaged primes, $\bar p_n = {1\over n} \sum_{i=1}^n p_i$, both the broad outline and the technical details of the arguments given and proofs presented are quite radically distinct.

\bigskip

\bigskip
\noindent
{\sc Date:} 19 July 2025; \LaTeX-ed \today

Oppermann; Brocard; Firoozbakht; Fourges; Nicholson; Farhadian.

}

\maketitle
\parindent0pt
\parskip7pt
\def\Mgap{{{\hbox{max-$\theta$-gap}}}}
\def\mgap{{{\hbox{min-$\theta$-gap}}}}
\def\O{{\mathcal{O}}}

\clearpage
\section{Introduction and background}

Consider the well-known sequence defined by evaluating the first Chebyshev function at the  $n^{th}$ prime:
\begin{equation}
\theta_n = \theta(p_n) =  \sum_{i=1}^n \ln p_i = \ln\left([p_n]\#\right).
\end{equation}
That is
\begin{equation}
\theta_n \in \{\ln 2, \ln 6, \ln 30, \ln 210, \ln 2310, \dots\} \quad\hbox{for}\quad n\in\{1,2,3,4,5, \dots\}.
\end{equation}
Quite a lot is already known about this sequence~\cite{Dusart:1998-thesis, Massias+Robin, Dusart:p_n}.
For instance the asymptotic behaviour is known to be $\theta_n \sim p_n \sim n \ln n$.
More precisely (Dusart~\cite{Dusart:1998-thesis}):
\begin{equation}
\theta_n =
n \left\{\ln n+ \ln\ln n - 1 + {\ln\ln n-2\over\ln n} 
+ \O\left(\ln\ln n\over(\ln n)^2\right) \right\}.
\end{equation}
Higher-order terms in the asymptotic expansion are also known
but are not of direct relevance to the present discussion.
Various other explicit upper and lower bounds are also known 
(see Dusart~\cite{Dusart:1998-thesis} and Massias+Robin~\cite{Massias+Robin}):
\begin{equation}
\theta_n \leq
n \left\{\ln n+ \ln\ln n - 1 + {\ln\ln n-2\over\ln n} 
\right\};\qquad
(n\geq 198).
\end{equation}
\begin{equation}
\theta_n \geq n \;[\ln (n\ln n) - 1] ;\qquad\qquad\qquad
(n\geq 5017).
\end{equation}
\begin{equation}
\theta_n \leq
n \; \ln (n\ln n);\qquad\qquad\qquad
(n\geq 3).
\end{equation}
\begin{equation}
\theta_n \geq
n \; \ln n;\qquad\qquad\qquad\qquad
(n\geq 13).
\end{equation}

Since $\theta_n \sim p_n$ these two sequences are very closely correlated.
Furthermore, since the relation, $p_n \longleftrightarrow \theta_n$ is easily invertible. $p_n = \exp(\theta_n - \theta_{n-1})$, these two sequences must ultimately encode identical information. 
Nevertheless, as discussed below, the two sequences also exhibit significant differences.
Specifically, I shall demonstrate below that the  sequence $\theta_n$ \emph{satisfies} suitably formulated $\theta$-analogues of the usual Cramer, Andrica, Legendre, Oppermann, Brocard, Firoozbakht, Fourges, Nicholson, and  Fahadian
conjectures. (So these $\theta$-analogues are not conjectures, they are instead theorems.) 
The crucial observation enabling this pleasant behaviour is the fact that the $\theta$-gaps $\mathfrak{g}_n = \theta_{n+1}- \theta_n = \ln p_{n+1}$ are in a suitable sense relatively small. Superficially these results bear  close resemblance to some recently derived results for the averaged primes~\cite{Visser:average-primes}, $\bar p_n = {1\over n} \sum_{i=1}^n p_i$, however both the broad outline and the technical details of the proofs are quite radically distinct.

\clearpage
\section{Standard and easy results}

In this section I will first introduce  a few basic tools:
\begin{eqnarray}
&p_n > n\ln n;  \qquad &(n\geq1; \hbox{ Rosser; \cite{Rosser1}});
\\
&p_n < n\ln(n\ln n);  \qquad\qquad &(n\geq6; \hbox{ Rosser+Schoenfeld; \cite{Rosser+Schoenfeld}});
\\
&p_n > n[\ln(n\ln n)-1];  \qquad\qquad &(n\geq2; \hbox{ Dusart; \cite{Dusart:p_n}}),
\\
&p_{n+1} < n[\ln(n\ln n)+1];  \qquad\qquad &(n\geq6; \hbox{ Visser; \cite{Visser:average-primes}});
\end{eqnarray}
Rosser and Schoenfeld~\cite{Rosser+Schoenfeld}  give $|\theta(x)-x| <{1\over2}\; {x\over\ln x} $ for $x\geq 563=p_{103}$. 
Checking small values of $x$ by explicit computation, $|\theta(x)-x| <{1\over2}\; \pi(x) $ for $x\geq 347=p_{69}$. 
But then, checking small values of $n$ by explicit computation, $|\theta_n - p_n | = |\theta(p_n)-p_n|<{1\over2} n$ for $n\geq 26$.
Thence, checking small integers by explicit computation one finds several simple explicit  bounds on the $\theta_n$:
\begin{eqnarray}
&\theta_n > n\left[\ln n-{1\over2}\right];  \qquad &(n\geq1);
\\
&\theta_n < n\left[\ln(n\ln n)+{1\over2}\right];  \qquad\qquad &(n\geq2);
\\
&\theta_n > n\left[\ln(n\ln n)-{3\over2}\right];  \quad\qquad\qquad &(n\geq2).
\end{eqnarray}

From Rosser and Schoenfeld~\cite{Rosser+Schoenfeld} we also have the asymmetric bounds
\begin{equation}
(x\geq 563=p_{103}) \qquad  
x\left(1- {1\over2\ln x}\right) < \theta(x) < x\left(1+ {1\over2\ln x}\right) 
\qquad (x>1).
\end{equation}
This implies
\begin{equation}
(n\geq {103}=\pi(563)) \qquad  
p_n\left(1- {1\over2\ln p_n}\right) < \theta_n < p_n\left(1+ {1\over2\ln p_n}\right) 
\qquad (n\geq 1).
\end{equation}
Many (much) tighter bounds along these lines (and similar lines) are known~\cite{Dusart:2018, Dusart:2010, Axler:2018, Axler:2019, little-book, unsolved, exponentially, dlVP, exponential-bounds} --- but these particular bounds will be sufficient for current purposes.

\clearpage
\section{The $\theta$-gaps:  $\mathfrak{g}_n =  \theta_{n+1}- \theta_n = \ln p_{n+1}$}
Consider the gaps between the $\theta_n$, specifically 
\begin{equation}
\mathfrak{g}_n =  \theta_{n+1}- \theta_n = \ln p_{n+1}.
\end{equation}
Note that both the $\theta_n$ themselves, and the gaps $\mathfrak{g}_n$, are monotone increasing, so this sequence is even smoother than the averaged primes $\bar p_n = {1\over n}\sum_{i=1}^n p_i$, for which the gaps in the averages, $\bar g_n = \bar p_{n+1}-\bar p_n$, are not monotone~\cite{Visser:average-primes}.

We find it useful to consider the relative gaps $\mathfrak{g}_n /\theta_n$. 
First note
\begin{equation}
{\mathfrak{g}_n \over \theta_n} = {\ln p_{n+1}  \over \theta_n}< {\ln p_{n+1}  \over p_n(1-{1\over 2 \ln p_n})}
< {\ln p_n\over p_n} \left( \ln p_{n+1} \over \ln p_n - {1\over2} \right);
\qquad (n \geq 103).
\end{equation}
Then apply the elementary Bertrand--Chebyshev result $p_{n+1}<2p_n$ to get 
\begin{equation}
{\mathfrak{g}_n \over \theta_n} < {\ln p_n\over p_n} \left\{ \ln p_{n} + \ln 2\over \ln p_n - {1\over2} \right\};
\qquad (n \geq 103).
\end{equation}
The quantity in braces is certainly less than 2 for $p_n \geq7$, that is for $n\geq 4$. 

So we certainly have
\begin{equation}
{\mathfrak{g}_n \over \theta_n} < 2 \;{\ln p_n\over p_n};
\qquad (n \geq 103).
\end{equation}
Explicitly checking smaller integers we find
\begin{equation}
{\mathfrak{g}_n \over \theta_n} < 2\; {\ln p_n\over p_n}; \qquad (n \geq 3).
\end{equation}

Working in the other direction
\begin{equation}
{\mathfrak{g}_n \over \theta_n} = {\ln p_{n+1}  \over \theta_n}> {\ln (p_{n}+2)  \over p_n(1+{1\over 2 \ln p_n})}
> {\ln p_n\over p_n} \left\{ \ln (p_{n} +2) \over \ln p_n + {1\over2} \right\};
\qquad (n \geq 1).
\end{equation}
The quantity in braces is easily seen to have an absolute minimum at $n=7$ where it takes on the value ${\ln19\over \ln 17+{1\over2}} = 0.8833634920$. So we certainly have
\begin{equation}
{\mathfrak{g}_n \over \theta_n}  > {22\over25} \; {\ln p_n\over p_n}; \qquad (n \geq 1).
\end{equation}

It is sometimes preferable to work directly in terms of $n$ rather than $p_n$ and write
\begin{equation}
{\mathfrak{g}_n \over \theta_n} = {\ln p_{n+1}  \over \theta_n}
< { \ln( n[\ln(n\ln n)+1]) \over n[\ln n-{1\over2}]}
= {1\over n} \left\{ \ln( n[\ln(n\ln n)+1]) \over \ln n-{1\over2}\right\}
\qquad (n \geq 6).
\end{equation}
The quantity in braces is certainly less than $1.8={9\over5}$ for $n\geq 17$, so
\begin{equation}
{\mathfrak{g}_n \over \theta_n} < {9/5\over n}; \qquad (n \geq 17).
\end{equation}
Explicitly checking smaller integers one finds
\begin{equation}
{\mathfrak{g}_n \over \theta_n} < {9/5\over n}; \qquad (n \geq 1).
\end{equation}
Working in the other direction
\begin{equation}
{\mathfrak{g}_n \over \theta_n} = {\ln p_{n+1}  \over \theta_n}
> {\ln( p_{n}  +2) \over \theta_n} > {\ln(n\ln n +2 ) \over n\left[\ln(n\ln n)+{1\over2}\right]}
= {1\over n} \; \left\{\ln(n\ln n+2) \over \ln(n\ln n)+{1\over2}\right\}
\qquad (n \geq 2).
\end{equation}
The quantity in braces is certainly greater than $0.9={9\over10}$ for $n\geq 22$, so
\begin{equation}
{\mathfrak{g}_n \over \theta_n} > {9/10\over n}; \qquad (n \geq 22).
\end{equation}
Explicitly checking smaller integers
\begin{equation}
{\mathfrak{g}_n \over \theta_n} > {9/10\over n}; \qquad (n \geq 1).
\end{equation}

Collecting results, we have relatively tight bounds on the relative $\theta$ gaps
\begin{equation}
(n\geq 1) \qquad
{22\over25} \; {\ln p_n\over p_n}< {\mathfrak{g}_n \over \theta_n} < 2\; {\ln p_n\over p_n}; 
\qquad (n \geq 3).
\end{equation}
\begin{equation}
(n \geq 1) \qquad 
{9/10\over n} < {\mathfrak{g}_n \over \theta_n} < {9/5\over n}; 
\qquad (n \geq 1).
\end{equation}
We will re-use these bounds many times in the discussion below.

\clearpage
\section{Counting the $\theta_n$}
\def\pit{{\pi_\theta}}
Define
\begin{equation}
\pit(x):= \#\{ i: \theta_i \leq x\}.
\end{equation}
This function counts the $\theta_n$.  Equivalently
\begin{equation}
\pit(x) = \sum_{i=1}^\infty \Heaviside(x-\theta_i); \qquad\qquad \Heaviside(0)=1.
\end{equation}
Let us now define
\begin{equation}
M_* = \max_{i\geq 1} \left(\theta_i\over p_i\right); \qquad
m_* = \min_{i\geq 1} \left(\theta_i\over p_i\right).
\end{equation}
That these two quantities exist is already guaranteed by the relatively simple bounds
(Rosser and Schoenfeld~\cite{Rosser+Schoenfeld}):
\begin{equation}
(n\geq {103}=\pi(563)) \qquad  
\left(1- {1\over2\ln p_n}\right) < {\theta_n\over p_n}  < \left(1+ {1\over2\ln p_n}\right) 
\qquad (n\geq 1).
\end{equation}
We also know that all of the numbers $\theta_n/p_n$ are transcendental, and that these extrema are unique. 
In fact if $\theta_n/p_n= \theta_m/p_m$ then $m=n$. To see this, first rearrange into the form $p_m \theta_n= p_n\theta_m$, and then exponentiate: $(\prod_{i=1}^n p_i)^{p_m} = (\prod_{i=1}^m p_i)^{p_n}$. But then, unless $n=m$ the LHS and RHS  will contain different collections of prime factors, which is impossible. 

Now we additionally  know $\theta_i/p_i < 1$ at least up to $p_i < 1.39\times10^{17} $, 
and that globally $M_* < 1 + 7.5\times10^{-7}$~(Platt+Trudgian~\cite{sign-change}). As for $m_*$, we have the relatively boring, (because it is rather weak), result that $m_* = \theta_1/p_1 = {\ln 2\over2} = 0.3465735903$. (We shall soon see how to improve on this.) For now, note that $ m_* p_i \leq \theta_i \leq M_* p_i$. 
Thence we see $\{ i: \theta_i \leq x\} \subseteq \{ i: m_* p_i \leq x\} = \{ i: p_i \leq x/m_*\}$, and so:
\begin{equation}
\pit(x) = \#\{ i: \theta_i \leq x\} \leq \#\{ i: p_i \leq x/m_*\} = \pi(x/m_*).
\end{equation}
Similarly, we see $\{ i: \theta_i \leq x\} \supseteq \{ i: M_* p_i \leq x\} = \{ i: p_i \leq x/M_*\}$, and so:
\begin{equation}
\pit(x) = \#\{ i: \theta_i \leq x\} \geq \#\{ i: p_i \leq x/M_*\} = \pi(x/M_*).
\end{equation}
Collecting results, for all $x\geq 1$ we have:
\begin{equation}
\pi(x/M_*) \leq \pit(x) \leq \pi(x/m_*).
\end{equation}
While the lower bound is reasonably tight, the upper bound is quite relaxed. 

Alternatively:
\begin{equation}
\pit(xm_*) \leq \pi(x) \leq \pit(xM_*).
\end{equation}
Note that the upper bound is now reasonably tight, whereas the lower bound is now quite relaxed. 

Let us see how to improve these results. 
First, let us take $x > y\geq 1$, and note
\begin{equation}
\pit(x)-\pit(y)= \#\{ \theta_i:  y < \theta_i \leq x \} > {x-y\over \Mgap{(x,y)}},
\end{equation}
where $ \Mgap(x,y)$ is the maximum $\theta$-gap that is completely contained within the interval $(y,x)$. Explicitly, noting that $p_i\leq \theta_i/m_*$, we have
\begin{equation}
\Mgap(x,y) = \max_{y<\theta_i< x} \left\{\ln p_i\right\} <  \max_{y<\theta_i< x} \left\{\ln[\theta_i/m_*]\right\}
= \max_{y<\theta_i< x} \left\{\ln\theta_i \right\} - \ln m_*.
\end{equation}
Finally, note that\; $\max_{y<\theta_i< x} \left\{\ln\theta_i \right\} < \ln x$, to obtain
\begin{equation}
\Mgap(x,y) < \ln x - \ln m_* =\ln x + \ln(2/\ln 2),
\end{equation}
and consequently we have the relatively weak lower bound that
\begin{equation}
\pit(x)-\pit(y) >  {x-y\over \ln x + \ln(2/\ln 2)}.
\end{equation}

Similarly to develop an upper bound we note that 
\begin{equation}
\pit(x)-\pit(y)= \#\{ \theta_i:  y < \theta_i \leq x \}  < {x-y\over \mgap{(x,y)}},
\end{equation}
where $ \mgap(x,y)$ is the minimum $\theta$-gap that is completely contained within the interval $(y,x)$. Explicitly, noting that $p_i\geq \theta_i/M_*$, we have
\begin{equation}
\mgap(x,y) = \min_{y<\theta_i< x} \left\{\ln p_i\right\} >  \min_{y<\theta_i< x} \left\{\ln[\theta_i/M_*]\right\}
= \min_{y<\theta_i< x} \left\{\ln\theta_i \right\} - \ln M_*.
\end{equation}
Finally, note $\min_{y<\theta_i< x} \left\{\ln\theta_i \right\} > \ln y$, to obtain
\begin{equation}
\mgap(x,y) > \ln y - \ln M_*,
\end{equation}
and consequently 
\begin{equation}
 \pit(x)-\pit(y) <  {x-y\over \ln y -\ln M_*}.
\end{equation}

Collecting results we have the practical and efficient bounds
\begin{equation}
 {x-y\over \ln x -\ln m_*} < \pit(x)-\pit(y) <  {x-y\over \ln y -\ln M_*},
\end{equation}
with
\begin{equation}
\ln m_* = \ln(\half\ln 2) = - 1.059660101; \qquad \ln M_* < \ln(1 + 7.5\times10^{-7}) <  7.5\times10^{-7}.
\end{equation}
Explicitly for $x> y\geq 1$ we have
\begin{equation}
 {x-y\over \ln x +1.059660101} < \pit(x)-\pit(y) <  {x-y\over \ln y - 7.5\times10^{-7}}.
\end{equation}
We shall have occasion to invoke these bounds multiple times in the discussion below.
(These bounds are almost certainly not the best that could in principle be developed, but they will suffice for our intended purposes.)

\clearpage
\section{$\theta$-analogues of some standard conjectures}

Given the rather tight bound on the $\theta$ gaps discussed above, it is now in very many cases relatively easy to formulate provable $\theta$-analogues of various standard conjectures. (Often the only major difficulty lies in designing and formulating a suitable $\theta$-analogue, and specifying a suitable domain of validity.) 
Below I shall present  a few examples.

\subsection{$\theta$-analogue of Cramer}

The ordinary Cramer conjecture is the hypothesis that the ordinary prime gaps satisfy $p_{n+1}-p_n = \O([\ln p_n]^2)$. 
The $\theta$-analogue would be $\theta_{n+1}-\theta_n = \O([\ln \bar p_n]^2) $.
But since we have as an identity that $\theta_{n+1}-\theta_n = \ln p_n $  this is a triviality.

\subsection{$\theta$-analogue of Andrica}

The ordinary Andrica conjecture is the hypothesis that the primes $p_n$ satisfy
\begin{equation}
\sqrt{p_{n+1}}- \sqrt{p_n} < 1; \qquad (n \geq 1),
\end{equation}
The $\theta$-analogue would be that the $\theta_n$ satisfy $\sqrt{\theta_{n+1}}- \sqrt{\theta_n} < 1$. But this is easily checked to be true and in fact much more can be said.

Using our previous results we compute
\begin{equation}
\sqrt{\theta_{n+1}} - \sqrt{\theta_n} = { \theta_{n+1}- \theta_n\over \sqrt{\theta_{n+1}} + \sqrt{\theta_n} } 
< { \ln p_{n+1} \over 2 \sqrt{\theta_n} } < { \ln(2 p_n) \over 2 \sqrt{ p_n(1-{1\over2\ln p_n}) } }.
\end{equation}
That is
\begin{equation}
\sqrt{\theta_{n+1}} - \sqrt{\theta_n} 
< {\ln p_n\over \sqrt{p_n}} \left\{ \sqrt{\ln p_n} (1+{\ln 2\over\ln p_n}) \over 2\sqrt{\ln p_n-{1\over2}}\right\}.
\end{equation}
But  for $n \geq 3$,  the quantity in braces lies in the range $(0,1)$ and so in this range we certainly have
\begin{equation}
\sqrt{\theta_{n+1}} - \sqrt{\theta_n}  <  {\ln p_n\over \sqrt{p_n}} ; \qquad (n\geq 3).
\end{equation}
Explicitly checking smaller integers
\begin{equation}
\sqrt{\theta_{n+1}} - \sqrt{\theta_n}  <  {\ln p_n\over \sqrt{p_n}} ; \qquad (n\geq 2).
\end{equation}

\clearpage
Note this is asymptotically much stronger than just a constant bound, and we could replace it with  a considerably weaker statement closer in spirit to the usual Andrica conjecture:
\begin{equation}
\sqrt{\theta_{n+1}} - \sqrt{\theta_n}  \leq  \sqrt{\theta_{2}} - \sqrt{\theta _1} = \sqrt{\ln 6} -\sqrt{\ln 2} = {0.5060115878}; \qquad (n\geq 1).
\end{equation}

Thence in particular
\begin{equation}
\sqrt{\theta_{n+1}} - \sqrt{\theta_n}  \leq 1; \qquad (n\geq 1).
\end{equation}
So the $\theta$-analogue of Andrica is unassailably true.\footnote{Effectively one has ``changed the problem'' by replacing a hard problem with something tractable. For other examples along these lines see for instance reference~\cite{Andrica-variants}.}

\subsection{$\theta$-analogue of Legendre}

The ordinary Legendre conjecture is the hypothesis that the ordinary primes $p_n$ satisfy $\pi([m+1]^2)> \pi(m^2)$. 
The $\theta$-analogue of Legendre would be that $\theta_n$ satisfy $\pit([m+1]^2)> \pit(m^2)$. But this is easily checked to be true and in fact much more can be said.

Applying our lower bound on $\pit(x)-\pit(y)$ we see that for $n\geq 1$ we have
\begin{equation}
\pit([n+1]^2)-\pit(n^2) > {[n+1]^2- n^2 \over \ln([n+1]^2) +1.059660101}
=
{2n+1\over 2 \ln(n+1)+1.059660101}.
\end{equation}
But this is manifestly greater than unity for $n\geq 1$.
So the $\theta$-analogue of Legendre  is manifestly true.
(In fact for $n\gg1$ there will be \emph{many}  $\theta_n$  between consecutive squares.) 

\subsection{$\theta$-analogue of Oppermann}

The ordinary Oppermann conjecture is the hypothesis that the ordinary primes $p_n$ satisfy $\pi(m[m+1])> \pi(m^2)> \pi(m[m-1])$ for integer $m\geq 2$. This is completely equivalent to demanding that $\pi([m+{1\over2}]^2)> \pi(m^2)> \pi([m-{1\over2}]^2)$ for $m\geq 2$. 
 
Then the  most compelling $\theta$-analogue of Oppermann would be that the $\theta_n$ satisfy 
$\pit([m+{1\over2}]^2)> \pit(m^2) > \pit([m-{1\over2}]^2)$. But (adapting the argument presented above for Legendre) this is easily checked to be true, and in fact much more can be said.

Note that for $m>1$ we have the lower bound
\begin{equation}
\pit([m+\half]^2)- \pit(m^2)  >{[m+\half]^2-m^2\over \ln([m+\half]^2)+1.059660101}
= {m+\quarter\over 2\ln(m+{1\over2}) +1.059660101}.
\end{equation}
But this is manifestly greater than unity for $m\geq 4$.
In fact, the observation that this is manifestly positive over the range of applicability ($m\geq 1$) is sufficient.

Similarly for $m\geq 2$ we have
\begin{equation}
\pit(m^2)  -\pit([m-\half]^2) > {m^2-[m-\half]^2\over \ln(m^2)+1.059660101}
= {m-\quarter\over 2\ln m +1.059660101}.
\end{equation}
But this is manifestly greater than unity for $m\geq 5$.
In fact, the observation that this is manifestly positive over the range of applicability ($m\geq 2$) is sufficient.  

In summary
\begin{equation}
(m\geq 1) \qquad \pit([m+\half]^2) > \pit(m^2)  >\pit([m-\half]^2); \qquad (m\geq 2).
\end{equation}
So the $\theta$-analogue of Oppermann  is manifestly true.
(In fact for $m\gg1$ there will be \emph{many}  $\theta_n$  between consecutive half-squares.) 

\subsection{$\theta$-analogue of Brocard}

The ordinary Brocard conjecture is the hypothesis that $\pi(p_{n+1}^2 )- \pi(p_n^2)\geq  4$ for $n\geq 2$. The ordinary Brocard conjecture is implied by the ordinary Oppermann conjecture.  Note that $p_{n+1}\geq p_n+2$, and that the ordinary Oppermann conjecture would imply the existence of at least one prime in each of the four open intervals:
\begin{equation}
\textstyle{
(p_n^2, [p_n+{1\over2}]^2),  \;\; ([p_n+{1\over2}]^2, [p_n+1]^2),  \;\; 
([p_n+1]^2, [p_n+{3\over2}]^2),  \;\;  ([p_n+{3\over2}]^2, [p_n+2]^2).
}
\end{equation}
Generalizing this, and checking the case $n=1$ by hand, the ordinary Oppermann conjecture would imply a generalization of the ordinary Brocard conjecture: 
\begin{equation}
\pi(p_{n+1}^2 )- \pi(p_n^2)\geq  2 g_n; \qquad (n \geq 1).
\end{equation}

Then one compelling $\theta$-analogue of Brocard would be to demand that the $\theta_n$ satisfy $\pit(p_{n+1}^2)- \pit(p_n^2) \geq 4$, or more generally $\pit(p_{n+1}^2)- \pit(p_n^2) \geq 2 g_n$, for suitable values of $n$. But (adapting the argument presented above for Legendre) this is easily checked to be true, and in fact much more can be said.

Observe
\begin{equation}
\pit(p_{n+1}^2)- \pit(p_{n}^2)  >{p_{n+1}^2- p_{n}^2\over \ln(p_{n+1}^2)+1.059660101}
= {(p_{n+1}+p_{n})(p_{n+1}-p_{n})\over 2\ln(p_{n+1}) +1.059660101}
\end{equation}
\begin{equation}
= g_n \left\{ p_{n+1}+p_{n}\over 2\ln(p_{n+1}) +1.059660101\right\}
\end{equation}
The quantity in braces is greater than 2 for $p_n\geq 5 = p_3$, that is, $n \geq 3$. That is
\begin{equation}
\pit(p_{n+1}^2 )- \pit(p_n^2)\geq  2 g_n; \qquad (n \geq 3).
\end{equation}
Checking smaller integers by direct computation
\begin{equation}
\pit(p_3^2)-\pit(p_2^2) = \pit(25)-\pit(9) = 10-5=5;\qquad g_2 = 5-3=2;
\end{equation}
\begin{equation}
\pit(p_2^2)-\pit(p_1^2) = \pit(9)-\pit(4) = 5-3=2; \qquad g_1=3-2=1.
\end{equation}
So we actually have
\begin{equation}
\pit(p_{n+1}^2 )- \pit(p_n^2)\geq  2 g_n; \qquad (n \geq 1).
\end{equation}
So a suitable $\theta$-analogue of Brocard is unassailably true.

\subsection{$\theta$-analogue of Firoozbakht}

The ordinary Firoozbakht conjecture is the hypothesis that the quantity $(p_n)^{1/n}$ is monotone decreasing for the ordinary primes $p_n$. The $\theta$-analogue of  Firoozbakht would be the hypothesis that $(\theta_n)^{1/n}$ is monotone decreasing for the $\theta_n$. But this is easily checked to be true (over a suitable domain to be determined below) and in fact much more can be said.

Note
\begin{equation}
(\theta_{n+1})^{1/(n+1)} < (\theta_n)^{1/n} \quad \iff \quad  n \ln \theta_{n+1} < (n+1) \ln \theta_n
\end{equation}
So let us compute
\begin{equation}
Q:= n \ln \theta_{n+1}-(n+1)\ln \theta_n = n \ln (\theta_{n} +\ln p_{n+1})-(n+1)\ln \theta_n
\end{equation}
\begin{equation}
= n \ln(1+\ln p_{n+1}/\theta_n) - \ln \theta_n. 
\end{equation}

This quantity is certainly less than
\begin{equation}
Q < {n \ln p_{n+1}\over  \theta_n} - \ln \theta_n 
< {n \ln (n[\ln(n\ln n)+1])  \over n[\ln(n\ln n) -{3\over2}]}- \ln (n\left[\ln(n\ln n) -{3\over2}\right]) 
\end{equation}
\begin{equation}
={\ln (n[\ln(n\ln n)+1])  \over \ln(n\ln n) -{3\over2}}- \ln (n\left[\ln(n\ln n) -{3\over2}\right]).
\end{equation}
This last bound becomes negative for $n \geq 9$. Checking smaller integers by direct computation shows $Q<0$ for $n\geq 4$, so the $\theta$-analogue of Firoozbakht is unassailably true for all $n\geq 4$.  (Indeed the sequence $[\theta_n]^{1/n}$ has a global maximum at $n=4$, the first few terms in that sequence are anomalous.)

\subsection{$\theta$-analogue of Fourges}

The ordinary Fourges conjecture is typically presented in terms of first rearranging the ordinary Firoozbakht conjecture into the form
\begin{equation}
\left( \ln p_{n+1} \over \ln p_n \right) < \left(1+{1\over n}\right);\qquad (n \geq 1),
\end{equation}
and then \emph{weakening} it by making the \emph{less restrictive} demand that
\begin{equation}
\left( \ln p_{n+1} \over \ln p_n \right) <  e^{1/n};\qquad (n \geq 1),
\end{equation}

To obtain a $\theta$-analogue of the Fourges conjecture we first rearrange the $\theta$-analogue Firoozbakht conjecture into the form
\begin{equation}
\left( \ln \theta_{n+1} \over \ln \theta_n \right) < \left(1+{1\over n}\right);\qquad (n \geq 4),
\end{equation}
and then \emph{weaken} it by making the \emph{less restrictive} demand that
\begin{equation}
\left( \ln \theta_{n+1} \over \ln \theta_n \right) <  e^{1/n};\qquad (n \geq 4).
\end{equation}
Since this is a weakening, the $\theta$-analogue of Fourges is unassailably true for all $n\geq 4$.  Checking smaller integers by explicit computation, the $\theta$-analogue of Fourges is unassailably true for all $n\in \{1\} \cup\{3,\infty\}$. That is,  $\theta$-Fourges is true for all positive integers except 2. 

\clearpage
\subsection{$\theta$-analogue of Nicholson}

The ordinary Nicholson conjecture is typically presented in terms of first rearranging the ordinary Firoozbakht conjecture into the form 
\begin{equation}
\left( p_{n+1} \over p_n \right)^n <  p_n;\qquad (n \geq 1),
\end{equation}
and then (slightly) \emph{strengthening} it by making the \emph{more restrictive} demand that
\begin{equation}
\left( p_{n+1} \over p_n \right)^n <  n \ln n;\qquad (n \geq 1).
\end{equation}

To obtain a $\theta$-analogue of the Nicholson conjecture we first rearrange 
the $\theta$-analogue Firoozbakht conjecture into the form
\begin{equation}
\left( \theta_{n+1} \over \theta_n \right)^n <  \theta_n;\qquad (n \geq 4),
\end{equation}
and then, noting that for $n\geq 1$ we have $\theta_n > n[\ln n -\half]$,  (slightly) \emph{strengthen} it by making the \emph{more  restrictive} demand that
\begin{equation}
\left( \theta_{n+1} \over \theta_n \right)^n <   n[\ln n -\half] ;\qquad (n \geq n_0),
\end{equation}
for some $n_0$ yet to be determined.
That is, our proposed form of the prime-average analogue of Nicholson is equivalent to 
\begin{equation}
n \ln \left( \theta_{n+1} \over \theta_n \right) <  \ln\left(  n[\ln n -\half]  \right) ;\qquad (n \geq n_0).
\end{equation}

Since
\begin{equation}
n \ln \left( \theta_{n+1} \over \theta_n \right) = n \ln \left(1+ {\mathfrak{g}_{n} \over \theta_n} \right)
<  n \ln \left(1+{9/5\over n}\right) < 9/5,
\end{equation}
with the inequalities holding for $n\geq 1$, we see that our proposed version of the prime-average analogue of the Nicholson conjecture will hold whenever
\begin{equation}
\ln\left(  n[\ln n -\half]  \right) > 9/5,
 \end{equation}
that is, whenever $n \geq 4$. Finally, explicitly checking all integers below 4 we have an explicit  prime-average analogue of Nicholson:
 \begin{equation}
\left( \theta_{n+1} \over \theta_n \right)^n <  n[\ln n -\half]   ;\qquad (n \geq 3).
\end{equation}
 
\subsection{$\theta$-analogue of Farhadian}

\enlargethispage{20pt}
The ordinary Farhadian conjecture is typically presented in terms of  (slightly) \emph{strengthening} the ordinary Nicholson conjecture by making the even \emph{more restrictive} demand that
\begin{equation}
\left( p_{n+1} \over p_n \right)^n <  {p_n \ln n\over\ln p_n};\qquad (n \geq 1).
\end{equation}

To obtain a $\theta$-analogue of the Farhadian conjecture we again  (slightly) \emph{strengthen} the $\theta$-analogue of Nicholson conjecture by making  \emph{some} (slightly) more restrictive demand.  

There is a potential infinity of stronger demands that one might make, but we shall try to keep close to the spirit of the original Farhadian conjecture by demanding
\begin{equation}
\left( \theta_{n+1} \over \theta_n \right)^n <  {p_n\over \ln p_n} \: [\ln n -\half] ;\qquad (n \geq n_0).
\end{equation}
Since $p_n/\ln p_n <n$ for $n\geq 4$ this certainly (in that region) strengthens the $\theta$-analogue of Nicholson. 
This would be equivalent to
\begin{equation}
n \ln\left( \theta_{n+1} \over \theta_n \right) <  \ln\left({p_n\over  \ln p_n} \: [\ln n -\half] \right);\qquad (n \geq n_0).
\end{equation}

As previously noted 
\begin{equation}
n \ln \left( \theta_{n+1} \over \theta_n \right) = n \ln \left(1+ {\mathfrak{g}_{n} \over \theta_n} \right)
<  n \ln \left(1+{9/5\over n}\right) < 9/5,
\end{equation}
with these inequalities holding for $n\geq 1$. 

This now suggests we  consider the inequality
\begin{equation}
\ln\left({p_n\over \ln p_n} \; [\ln n -\half] \right) > 9/5,
\end{equation}
which holds for $n\geq 6$. Explicitly checking small integers
\begin{equation}
n \ln\left( \theta_{n+1} \over \theta_n \right) <  \ln\left({p_n\over  \ln p_n} \: [\ln n -\half] \right);\qquad (n \geq 5).
\end{equation}
Equivalently 
\begin{equation}
\left( \theta_{n+1} \over \theta_n \right)^n <  {p_n\over \ln p_n} \: [\ln n -\half] ;\qquad (n \geq 5).
\end{equation}
As expected, because the $\theta$-analogue of Farhadian is slightly stronger than the $\theta$-analogue of Nicholson, it is valid only on a slightly smaller domain.

\section{Discussion}\label{S:discussion}

We have explicitly seen above that the sequence $\theta_n = \theta(p_n) = \sum_{i=1}^n \ln p_i$ satisfies suitably defined $\theta$-analogues of 
the usual Cramer, Andrica, Legendre, Oppermann, Brocard, and Firoozbakht, Fourges, Nicholson, and Farhadian  hypotheses (in suitably defined domains). The key step in all cases was the extremely tight bound on the $\theta$ gaps: $\mathfrak{g_n} = \theta_{n+1}-\theta_n  = \ln p_{n+1} $. 
The only (relatively minor) potential difficulty lies in the design and formulation of truly compelling $\theta$-analogues of these hypotheses.

On the other hand, what does this tell us about the ordinary primes $p_n$?\\
Not as much as one might hope. While we certainly have the identity
\begin{equation}
g_n = p_{n+1}-p_n = \exp(\mathfrak{g}_n) - \exp(\mathfrak{g}_{n-1}) ,
\end{equation}
unfortunately naive applications of this equality do not  lead to anything useful. So even though the numerical evidence is quite compelling~\cite{Verifying1,Verifying2}, significant progress on the usual Cramer, Andrica, Legendre, Oppermann, Brocard, and Firoozbakht, Fourges, Nicholson, and Farhadian  conjectures for ordinary primes would require considerably more subtle arguments. 

\bigskip

\hrule\hrule\hrule

\clearpage
\addtocontents{toc}{\bigskip\hrule}

\setcounter{secnumdepth}{0}
\section[\hspace{14pt}  References]{}
%
\vspace{-50pt}

\end{document}